\documentclass{commat}

\usepackage{graphicx}

\title[A regularity criterion in multiplier spaces to Navier-Stokes equations]{A regularity criterion in multiplier spaces to Navier-Stokes equations via the gradient of one velocity component}

\author{Ahmad M. Alghamdi, Sadek Gala and Maria Alessandra Ragusa}

\affiliation{%
    Ahmad M. Alghamdi --
        Department of Mathematical Sciences, Faculty of Applied Sciences, Umm Al-Qura University, Saudi Arabia. P.O. Box 14035, Makkah 21955.
    \email{amghamdi@uqu.edu.sa}
    Sadek Gala --
        Dipartimento di Matematica e Informatica, Universit\`{a} di Catania, Viale Andrea Doria, 6 95125 Catania - Italy.
    \email{sadek.gala@univ-mosta.dz, sgala793@gmail.com}
    Maria Alessandra Ragusa --
        Dipartimento di Matematica e Informatica, Universit\`{a} di Catania, Viale Andrea Doria, 6 95125 Catania - Italy.  RUDN University, 6 Miklukho - Maklay St, Moscow, 117198, Russia.
    \email{mariaalessandra@unict.it}
    }

\abstract{
    In this paper, we study regularity of weak solutions to the incompressible Navier-Stokes equations in $\mathbb{R}^{3}\times (0,T)$. The main goal is to establish the regularity criterion via the gradient of one velocity component in multiplier spaces.
    }



\VOLUME{31}
\YEAR{2023}
\NUMBER{1}
\firstpage{81}
\DOI{https://doi.org/10.46298/cm.10267}

\begin{paper}

\section{Introduction}

In this paper we consider the following Cauchy problemma for the
incompressible Navier-Stokes equations in $\mathbb{R}^{3}\times (0,T)$
\begin{equation}
\left\{
\begin{array}
{c}
\partial _{t}u+\left(u\cdot \nabla\right) u-\Delta u+\nabla \pi = 0, \\
\nabla \cdot u = 0, \\
u(x,0) = u_{0}(x),
\end{array}
\right. \label{eq1}
\end{equation}
where $u = (u_{1}(x,t),u_{2}(x,t),u_{3}(x,t))$ denotes the unknown velocity
vector and $\pi = \pi (x,t)$ denotes the hydrostatic pressure respectively.
While $u_{0}$ is the prescribed initial data for the velocity with
properties $\nabla \cdot u_{0} = 0$.

The global existence of smooth solutions for the 3D incompressible
Navier-Stokes equations is one of the most outstanding open problemmas in
fluid mechanics. Different criteria for regularity of the weak solutions
have been proposed and many interesting results were established (see, for
example,~\cite{Bei}, \cite{FO}, \cite{FJNZ}, \cite{He1}, \cite{HG}, \cite{GG}, \cite{KOT}, \cite{Ser}, \cite{ZC}, \cite{Z3}, \cite{Z4}, \cite{Z5}, \cite{ZG1}, \cite{ZP1} and references therein).

Recently, many authors became interested in the regularity criteria
involving only one velocity component, or its gradient, even though most of
which are not scaling invariant (see, for example,~\cite{CT0}, \cite{CT1}, \cite{CZ}, \cite{FQ1}, \cite{FQ2}, \cite{He}, \cite{JZ1}, \cite{PP}, \cite{Zhe} and the references cited therein). In particular,
Zhou~\cite{Z1} showed that the solution is regular if
\begin{equation}
u_{3}\in L^{p}(0,T;L^{q}(\mathbb{R}^{3}))\text{\ \ with \ }\frac{2}{p}+\frac{3}{q} = \frac{1}{2},\text{\ \ }6<q\leq \infty . \label{eq02}
\end{equation}
Later, Cao and Titi~\cite{CT1} obtained the regularity criterion
\begin{equation}
u_{3}\in L^{p}(0,T;L^{q}(\mathbb{R}^{3}))\text{\ \ with \ }\frac{2}{p}+\frac{3}{q} = \frac{2}{3}+\frac{2}{3q},\text{\ \ }q>\frac{7}{2}. \label{eq03}
\end{equation}
Motivated by the above work, Zhou and Pokorn\'{y}~\cite{ZP1} showed the
following regularity condition
\begin{equation}
u_{3}\in L^{p}(0,T;L^{q}(\mathbb{R}^{3}))\text{\ \ with \ }\frac{2}{p}+\frac{3}{q} = \frac{3}{4}+\frac{1}{2q},\text{\ \ }q>\frac{10}{3}, \label{eq2}
\end{equation}
while the limiting case $u_{3}\in L^{\infty}(0,T;L^{\frac{10}{3}}(\mathbb{R}^{3}))$was covered in~\cite{JZ1}. Inspired by the work~\cite{Bei}, we are
interested in criteria involving the gradient of one velocity component $\nabla u_{3}$. In fact, He~\cite{He} first verified the following regularity
result
\begin{equation}
\nabla u_{3}\in L^{p}(0,T;L^{q}(\mathbb{R}^{3}))\text{\ \ with \ }\frac{2}{p%
}+\frac{3}{q} = 1,\text{\ \ }3\leq q\leq \infty . \label{eq04}
\end{equation}
The above result was significantly improved by Pokorn\'{y}~\cite{Pok} and
Zhou~\cite{Z2} independently (see also~\cite{ZP2}). More precisely, it reads
as follows
\begin{equation}
\nabla u_{3}\in L^{p}(0,T;L^{q}(\mathbb{R}^{3}))\text{\ \ with \ }\frac{2}{p%
}+\frac{3}{q} = \frac{3}{2},\text{\ \ }2\leq q<\infty . \label{eq05}
\end{equation}
Very recently, Ye~\cite{Ye2} improves the previous work of Zhou and Pokorn\'{y}~\cite{ZP2} by using of a~new anisotropic Sobolev inequality, and
proved the following regularity criterion
\begin{equation}
\nabla u_{3}\in L^{\frac{16q}{15q-23}}(0,T;L^{q}(\mathbb{R}^{3}))
\quad \textup{ with } q \in [2,3]. \label{eq100}
\end{equation}
Note that
\[
\frac{2}{\frac{16q}{15q-23}} + \frac{3}{q}
= \frac{15q+1}{8q}
> \frac{23}{12},
\quad \textup{for any } 2\leq q < 3.
\]
Consequently, (\ref{eq100}) can be regarded as a~further improvement of~\cite{ZP2}. Moreover, the endpoint case $q = 3$recovers the
result of~\cite{ZY}. For some other interesting regularity criteria, we
refer the readers to~\cite{Ye1}, \cite{Ye2} and references therein.

The purpose of this work is to extend the regularity criterion of weak
solutions in terms of one gradient of velocity component to the multiplier
space which is larger than the Lebesgue space. The method is based on the
following interpolation inequality
\begin{equation*}
\left\Vert \varphi \right\Vert _{L^{\gamma}}\leq C\left\Vert \partial
_{3}\varphi \right\Vert _{L^{\mu}}^{\frac{1}{3}}\left\Vert \nabla
_{h}\varphi \right\Vert _{L^{\lambda}}^{\frac{2}{3}},
\end{equation*}
where $\mu ,\lambda $and $\gamma $satisfy
\begin{equation*}
1\leq \mu ,\text{\ \ }\lambda <+\infty ,\text{\ \ }\frac{1}{\mu}+\frac{2}{%
\lambda}>1\text{\ and \ }1+\frac{3}{\gamma} = \frac{1}{\mu}+\frac{2}{%
\lambda}.
\end{equation*}
The detailed proof of this inequality can be found in the appendix of Cao
and Wu~\cite{CW}.

In order to prove our theorem, let us recall \ the definition of weak
solutions.

\begin{definition}
Let $T>0, u_{0} \in L^{2}(\mathbb{R}^{3})$ with $\nabla \cdot u_{0} = 0$ in
the sense of distributions. A measurable function $u(x,t)$ is called a~weak
solution to the Navier-Stokes equations (\ref{eq1}) on $[0,T]$ if the
following conditions hold:

\begin{enumerate}
\item $u(x,t)\in L^{\infty}(0,T;L^{2}(\mathbb{R}^{3}))\cap L^{2}(0,T;H^{1}(\mathbb{R}^{3}))$;

\item system (\ref{eq1}) is satisfied in the sense of distributions;

\item the energy inequality, that is,
\begin{equation*}
\left\Vert u(\cdot ,t)\right\Vert _{L^{2}}^{2}+2\int_{0}^{t}\left\Vert
\nabla u(\tau)\right\Vert _{L^{2}}^{2}d\tau \leq \left\Vert
u_{0}\right\Vert _{L^{2}}^{2}.
\end{equation*}
\end{enumerate}
\end{definition}

By a~strong solution, we mean that a~weak solution $u$ of the Navier-Stokes
equations (\ref{eq1}) satisfies
\begin{equation*}
(u(x,t),\theta (x,t))\in L^{\infty}(0,T;H^{1}(\mathbb{R}^{3}))\cap
L^{2}(0,T;H^{2}(\mathbb{R}^{3})).
\end{equation*}
It is well known that the strong solution is regular and unique.

For $\alpha \in \mathbb{R}$, the Homogeneous Sobolev Space $\dot{H}^{\alpha
}(\mathbb{R}^{3})$is the space of tempered distributions $f$for which
\begin{equation*}
\left\Vert f\right\Vert _{\dot{H}^{\alpha}} = \sqrt{\int_{\mathbb{R}%
^{3}}\left\vert \xi \right\vert ^{2\alpha}\left\vert \widehat{f}(\xi
)\right\vert ^{2}d\xi}<+\infty .
\end{equation*}
For Homogeneous Sobolev Spaces, we refer to the book~\cite{BCD}. For
instance, the following basic interpolation inequality holds:

\begin{lemma}
For $0<\alpha \leq \beta$, the space $L^{2}\cap \dot{H}^{\beta}$ is a
subset of $\dot{H}^{\alpha}$, and we have
\begin{equation}
\left\Vert f\right\Vert _{\dot{H}^{\alpha}}\leq \left\Vert f\right\Vert
_{L^{2}}^{1-\frac{\alpha}{\beta}}\left\Vert f\right\Vert _{\dot{H}^{\beta
}}^{\frac{\alpha}{\beta}}. \label{eq25}
\end{equation}
\end{lemma}

\begin{proof}
This is a~particular case of~\cite{BCD}, Proposition 1.32.
\end{proof}

We say that a~function belongs to the multiplier spaces $\dot{X}_{1+\alpha
} := M(\dot{H}^{\alpha}(\mathbb{R}^{3})\rightarrow \dot{H}^{-1}(\mathbb{R}^{3}))$if it maps, by pointwise multiplication, $\dot{H}^{\alpha}$in $\dot{H}^{-1}:$
\begin{equation*}
\dot{X}_{1+\alpha} = \left\{f\in \mathcal{S}^{\prime}(\mathbb{R}^{3}):\left\Vert fg\right\Vert _{\dot{H}^{-1}}\leq \left\Vert g\right\Vert _{\dot{H}^{\alpha}}\right\}.
\end{equation*}
$\dot{H}^{\alpha}(\mathbb{R}^{3})$denotes the homogeneous Sobolev space.
The space $\dot{X}_{1+\alpha}$has been characterized in~\cite{Maz1}, \cite{Maz2}
(see also~\cite{GL}). Now our regularity criterion for system (\ref{eq1})
reads

\begin{theorem}\label{th1}Let $u_{0}\in L^{2}(\mathbb{R}^{3})$ with $\nabla \cdot u_{0} = 0$
in the sense of distributions. Assume that $u$ is a~weak solution to system (\ref{eq1}).\ If $\nabla u_{3}$ satisfies the following condition
\begin{equation}
\nabla u_{3}\in L^{\frac{8}{3-4\alpha}}(0,T;\dot{X}_{1+\alpha}(\mathbb{R}^{3})),\text{\ for some \ }0\leq \alpha <\frac{3}{4}, \label{eq16}
\end{equation}
then the solution $u$ is regular on $(0;T]$.
\end{theorem}

\begin{remark}
Since $L^{\frac{3}{1+\alpha}}(\mathbb{R}^{3})\subset \dot{X}_{1+\alpha}(\mathbb{R}^{3})$ (see e.g.~\cite{ZG1} for details), it is clear that our
result improves that in~\cite{Ye2} and extend the regularity criterion (\ref{eq2}) from Lebesgue space $L^{\alpha}$ to multiplier space $\dot{X}_{1+\alpha}$.
\end{remark}

Thanks to
\begin{equation*}
\left\Vert f\right\Vert _{BMO}\leq C\left\Vert \nabla f\right\Vert _{\dot{X}_{1}}
\end{equation*}
(see e.g.~\cite[Proposition~2]{G08}), where $BMO$ denotes the homogeneous
space of bounded mean oscillations, it is easy to deduce the following
regularity criterion.

\begin{corollary}
Let $u_{0}\in L^{2}(\mathbb{R}^{3})$ with $\nabla \cdot u_{0} = 0$ in the
sense of distributions. Assume that $u$ is a~weak solution to system (\ref{eq1}).\ If $u_{3}$ satisfies the following condition
\begin{equation}
u_{3}\in L^{\frac{8}{3}}(0,T;BMO(\mathbb{R}^{3})), \label{eq17}
\end{equation}
then the solution $u$ is regular on $(0;T]$.
\end{corollary}

\begin{remark}
Since $L^{\infty}(\mathbb{R}^{3})\hookrightarrow BMO(\mathbb{R}^{3})$, our
result recovers the limiting case $q = \infty $ in (\ref{eq2}), that is,
\begin{equation*}
u_{3}\in L^{\frac{8}{3}}(0,T;L^{\infty}(\mathbb{R}^{3})).
\end{equation*}
Consequently, (\ref{eq17}) can be regarded as a~further improvement of the
previous work~\cite{ZP1}.
\end{remark}

\section{Proof of main result.}

In this section, under the assumptions of the Theorem~\ref{th1}, we prove
our main result. Before proving our result, we recall the following
muliplicative Sobolev imbedding inequality in the whole space $\mathbb{R}^{3} $(see, for example~\cite{CT1}):
\begin{equation}
\left\Vert f\right\Vert _{L^{6}}\leq C\left\Vert \nabla _{h}f\right\Vert
_{L^{2}}^{\frac{2}{3}}\left\Vert \partial _{3}f\right\Vert _{L^{2}}^{\frac{1%
}{3}}, \label{eq8}
\end{equation}
where $\nabla _{h} = (\partial _{x_{1}},\partial _{x_{2}})$is the horizontal
gradient operator. We are now give the proof of our main theorem.

\begin{proof}
To prove our result, it suffices to show that for any fixed $T>T^{\ast}$,
there holds
\begin{equation*}
\underset{0\leq t\leq T^{\ast}}{\sup}\left\Vert \nabla u(t)\right\Vert
_{L^{2}}^{2}\leq C_{T},
\end{equation*}
where $T^{\ast}$, which denotes the maximal existence time of a~strong
solution and $C_{T}$ is an absolute constant which only depends on $T$ and $u_{0}$.

The method of our proof is based on two major parts. The first one
establishes the bounds of $\left\Vert \nabla _{h}u\right\Vert _{L^{2}}^{2}$,
while the second gives the bounds of the $H^{1}-$norm of velocity $u$ in
terms of the results of part one.

\textbf{Step I.} Taking the inner product of (\ref{eq1})$_{1}$ with $-\Delta
_{h}u$, we obtain after integrating by parts that
\begin{equation}
\frac{1}{2}\frac{d}{dt}\left\Vert \nabla _{h}u\right\Vert
_{L^{2}}^{2}+\left\Vert \nabla \nabla _{h}u\right\Vert _{L^{2}}^{2} = \int_{\mathbb{R}^{3}}(u\cdot \nabla)u\cdot \Delta _{h}udx = I. \label{eq5}
\end{equation}
where $\Delta _{h} = \partial _{x_{1}}^{2}+\partial _{x_{2}}^{2}$ is the
horizontal Laplacian. For the notational simplicity, we set
\begin{align*}
\mathcal{L}^{2}(t) &= \underset{\tau \in \lbrack \Gamma ,t]}{\sup}\left\Vert \nabla _{h}u(\tau)\right\Vert _{L^{2}}^{2}+\int_{\Gamma
}^{t}\left\Vert \nabla \nabla _{h}u(\tau)\right\Vert _{L^{2}}^{2}d\tau , \\
\mathcal{J}^{2}(t) &= \underset{\tau \in \lbrack \Gamma ,t]}{\sup}\left\Vert \nabla u(\tau)\right\Vert _{L^{2}}^{2}+\int_{\Gamma
}^{t}\left\Vert \Delta u(\tau)\right\Vert _{L^{2}}^{2}d\tau ,
\end{align*}
for $t\in \lbrack \Gamma ,T^{\ast})$. In view of (\ref{eq16}), we choose $\epsilon >0$ to be precisely determined subsequently and then select $\Gamma
<T^{\ast}$ sufficiently close to $T^{\ast}$ such that for all $\Gamma \leq
t<T^{\ast}$,
\begin{equation}
\int_{\Gamma}^{t}\left\Vert \nabla u(\tau)\right\Vert _{L^{2}}^{2}d\tau
\leq \epsilon \ll 1. \label{eq19}
\end{equation}
Integrating by parts and using the divergence-free condition, it follows
that
\begin{align*}
I &\leq \int_{\mathbb{R}^{3}}\left\vert \nabla u_{3}\right\vert \left\vert
\nabla u\right\vert \left\vert \nabla _{h}u\right\vert dx \\
&\leq \left\Vert \left\vert \nabla u_{3}\right\vert \left\vert \nabla
u\right\vert \right\Vert _{\dot{H}^{-1}}\left\Vert \nabla _{h}u\right\Vert _{\dot{H}^{1}} \\
&\leq C\left\Vert \nabla u_{3}\right\Vert _{\dot{X}_{1+\alpha}}\left\Vert
\nabla u\right\Vert _{\dot{H}^{\alpha}}\left\Vert \nabla \nabla
_{h}u\right\Vert _{L^{2}} \\
&\leq C\left\Vert \nabla u_{3}\right\Vert _{\dot{X}_{1+\alpha}}\left\Vert
\nabla u\right\Vert _{L^{2}}^{1-\alpha}\left\Vert \nabla ^{2}u\right\Vert
_{L^{2}}^{\alpha}\left\Vert \nabla \nabla _{h}u\right\Vert _{L^{2}} \\
&\leq C\left\Vert \nabla u_{3}\right\Vert _{\dot{X}_{1+\alpha
}}^{2}\left\Vert \nabla u\right\Vert _{L^{2}}^{2(1-\alpha)}\left\Vert
\Delta u\right\Vert _{L^{2}}^{2\alpha}+\frac{1}{2}\left\Vert \nabla \nabla
_{h}u\right\Vert _{L^{2}}^{2},
\end{align*}
by Young's inequality and (\ref{eq25}). Inserting the above estimate into (\ref{eq5}) and integrating with respect to time, we deduce for every $\tau
\in \lbrack \Gamma ,t]$:
\newline
\resizebox{\textwidth}{!}{
\parbox{1.36\textwidth}{
\begin{align*}
\underset{\tau \in \lbrack \Gamma ,t]}{\sup}&\left\Vert \nabla _{h}u(\tau
)\right\Vert _{L^{2}}^{2}+\int_{\Gamma}^{t}\left\Vert \nabla \nabla
_{h}u(\tau)\right\Vert _{L^{2}}^{2}d\tau \\
&\leq\left\Vert \nabla _{h}u(\Gamma)\right\Vert
_{L^{2}}^{2}+C\int_{\Gamma}^{t}\left\Vert \nabla u_{3}(\tau)\right\Vert _{\dot{X}_{1+\alpha}}^{2}\left\Vert \nabla u(\tau)\right\Vert
_{L^{2}}^{2(1-\alpha)}\left\Vert \Delta u(\tau)\right\Vert
_{L^{2}}^{2\alpha}d\tau \\
&\leq \left\Vert \nabla _{h}u(\Gamma)\right\Vert _{L^{2}}^{2}+C\left(
\underset{\tau \in \lbrack \Gamma ,t]}{\sup}\left\Vert \nabla u(\tau
)\right\Vert _{L^{2}}^{\frac{3}{2}-2\alpha}\right) \int_{\Gamma
}^{t}\left\Vert \nabla u_{3}(\tau)\right\Vert _{\dot{X}_{1+\alpha
}}^{2}\left\Vert \nabla u(\tau)\right\Vert _{L^{2}}^{\frac{1}{2}}\left\Vert
\Delta u(\tau)\right\Vert _{L^{2}}^{2\alpha}d\tau \\
&\leq C+C\left(\underset{\tau \in \lbrack \Gamma ,t]}{\sup}\left\Vert
\nabla u(\tau)\right\Vert _{L^{2}}^{\frac{3}{2}-2\alpha}\right) \left(
\int_{\Gamma}^{t}\left\Vert \nabla u_{3}(\tau)\right\Vert _{\dot{X}_{1+\alpha}}^{\frac{8}{3-4\alpha}}d\tau\right) ^{\frac{3}{4}-\alpha
}\left(\int_{\Gamma}^{t}\left\Vert \nabla u(\tau)\right\Vert
_{L^{2}}^{2}d\tau\right) ^{\frac{1}{4}}\left(\int_{\Gamma}^{t}\left\Vert
\Delta u(\tau)\right\Vert _{L^{2}}^{2}d\tau\right) ^{\alpha} \\
&\leq C+C\mathcal{J}^{\frac{3}{2}-2\alpha}(t)\left(\int_{\Gamma
}^{t}\left\Vert \nabla u_{3}(\tau)\right\Vert _{\dot{X}_{1+\alpha}}^{\frac{%
8}{3-4\alpha}}d\tau\right) ^{\frac{3}{4}-\alpha}\epsilon ^{\frac{1}{4}}\mathcal{J}^{2\alpha}(t) \\
&\leq C+C\epsilon ^{\frac{1}{4}}\mathcal{J}^{\frac{3}{2}}(t),
\end{align*}
}}
which leads to
\begin{equation}
\mathcal{L}^{2}(t)\leq C+C\epsilon ^{\frac{1}{4}}\mathcal{J}^{\frac{3}{2}}(t). \label{eq4}
\end{equation}

\textbf{Step II.} Now, we will establish the bounds of $H^{1}-$norm of the
velocity field. In order to do it, taking the inner product of (\ref{eq1})$_{1}$ with $-\Delta u$ in $L^{2}(\mathbb{R}^{3})$. Then, integration by
parts gives the following identity:
\begin{equation*}
\frac{1}{2}\frac{d}{dt}\left\Vert \nabla u\right\Vert
_{L^{2}}^{2}+\left\Vert \Delta u\right\Vert _{L^{2}}^{2} = \int_{\mathbb{R}^{3}}(u\cdot \nabla)u\cdot \Delta udx
\end{equation*}
Integrating by parts and using the divergence-free condition, one can easily
deduce that (see e.g.~\cite{ZP1})
\begin{align*}
\int_{\mathbb{R}^{3}}(u\cdot \nabla)u\cdot \Delta udx 
&\leq C\int_{\mathbb{%
R}^{3}}\left\vert \nabla _{h}u\right\vert \left\vert \nabla u\right\vert
^{2}dx\leq C\left\Vert \nabla _{h}u\right\Vert _{L^{2}}\left\Vert \nabla
u\right\Vert _{L^{4}}^{2} \\
&\leq C\left\Vert \nabla _{h}u\right\Vert _{L^{2}}\left\Vert \nabla
u\right\Vert _{L^{2}}^{\frac{1}{2}}\left\Vert \nabla u\right\Vert _{L^{6}}^{\frac{3}{2}} \\
&\leq C\left\Vert \nabla _{h}u\right\Vert _{L^{2}}\left\Vert \nabla
u\right\Vert _{L^{2}}^{\frac{1}{2}}\left\Vert \nabla _{h}\nabla u\right\Vert
_{L^{2}}\left\Vert \Delta u\right\Vert _{L^{2}}^{\frac{1}{2}}
\end{align*}
by H\"{o}lder's inequality, Nirenberg-Gagliardo's interpolation inequality
and (\ref{eq8}). Integrating this last inequality in time, we deduce that
for all $\tau \in \lbrack \Gamma ,t]$
\newline
\resizebox{\textwidth}{!}{
\parbox{1.25\textwidth}{
\begin{align}
\mathcal{J}^{2}(t)
\le{ }&{ } \left\Vert \nabla u(\Gamma)\right\Vert _{L^{2}}^{2} \notag \\
&{\ } + C\underset{\tau \in \lbrack \Gamma ,t]}{\sup}\left\Vert \nabla _{h}u(\tau)\right\Vert _{L^{2}}\left(\int_{\Gamma}^{t}\left\Vert \nabla u(\tau)\right\Vert _{L^{2}}^{2}d\tau\right) ^{\frac{1}{4}}\left( \int_{\Gamma}^{t}\left\Vert \nabla \nabla _{h}u(\tau)\right\Vert _{L^{2}}^{2}d\tau\right) ^{\frac{1}{2}}\left(\int_{\Gamma}^{t}\left\Vert \Delta u(\tau)\right\Vert _{L^{2}}^{2}d\tau\right) ^{\frac{1}{4}} \notag
\\
\le{ }&{ } \left\Vert \nabla u(\Gamma)\right\Vert _{L^{2}}^{2}+2C\mathcal{L}(t)\epsilon ^{\frac{1}{4}}\mathcal{L}(t)\mathcal{J}^{\frac{1}{2}}(t) \notag
\\
={ }&{ } \left\Vert \nabla u(\Gamma)\right\Vert _{L^{2}}^{2}+C\epsilon ^{\frac{1}{%
4}}\mathcal{L}^{2}(t)\mathcal{J}^{\frac{1}{2}}(t). \label{eq26}
\end{align}
}}

Inserting (\ref{eq4}) into (\ref{eq26}) and taking $\epsilon $ small enough,
then it is easy to see that for all $\Gamma \leq t<T^{\ast}$, there holds
\begin{equation*}
\mathcal{J}^{2}(t)\leq \left\Vert \nabla u(\Gamma)\right\Vert
_{L^{2}}^{2}+C\epsilon ^{\frac{1}{4}}\mathcal{J}^{\frac{1}{2}}(t)+C\epsilon
^{\frac{1}{2}}\mathcal{J}^{2}(t)<\infty ,
\end{equation*}
which proves
\begin{equation*}
\underset{\Gamma \leq t<T^{\ast}}{\sup}\left\Vert \nabla u(t)\right\Vert
_{L^{2}}^{2}<+\infty .
\end{equation*}
This completes the proof of Theorem~\ref{th1}.
\end{proof}

\subsection*{Acknowledgments}

Part of the work was carried out while the second author was long-term
visitor at University of Catania. The hospitality of Catania University is
graciously acknowledged. This research is partially supported by Piano della
Ricerca 2016-2018 - Linea di intervento 2: ``Metodi variazionali ed equazioni
differenziali''. The second author wish to thank the support of ``RUDN
University Program 5-100''. The authors thank the anonymous referee for
careful reading of the manuscript, many valuable comments and helprooful
suggestions for its improvement.

\EditInfo{November 23, 2019}{May 29, 2020}{Diana Barseghyan}

\end{paper}
\begin{references}
\refer{Paper}{BL}
\Rauthor{J. Bergh and J. Löfström}
\Rtitle{Interpolation Spaces}
\Rseries{Grundlehren der mathematischen Wissenschaften [Fundamental
          Principles of Mathematical Sciences]}
\Rvolume{223}
\Rpublisher{Springer Berlin, Heidelberg}
\Ryear{1976}
\Rpages{207}

\refer{Book}{BCD}
\Rauthor{Bahouri, Hajer and Chemin, Jean-Yves and Danchin, Rapha\"{e}l}
\Rtitle{Fourier analysis and nonlinear partial differential equations}
\Rseries{Grundlehren der mathematischen Wissenschaften [Fundamental
              Principles of Mathematical Sciences]}
\Rvolume{343}
\Rpublisher{Springer, Heidelberg}
\Ryear{2011}
\Rpages{xvi+523}

\refer{Paper}{Bei}
\Rauthor{Beir\~{a}o da Veiga, H.}
\Rtitle{A new regularity class for the {N}avier-{S}tokes equations in {${\bf R}^n$}}
\Rjournal{Chinese Ann. Math. Ser. B}
\Rvolume{16}
\Ryear{1995}
\Rnumber{4}
\Rpages{407-412}

\refer{Paper}{CT0}
\Rauthor{Cao, Chongsheng and Titi, Edriss S.}
\Rtitle{Global regularity criterion for the 3{D} {N}avier-{S}tokes
              equations involving one entry of the velocity gradient tensor}
\Rjournal{Arch. Ration. Mech. Anal.}
\Rvolume{202}
\Ryear{2011}
\Rnumber{3}
\Rpages{919-932}

\refer{Book}{CT1}
\Rauthor{Cao, Chongsheng and Titi, Edriss S.}
\Rtitle{Regularity criteria for the three dimensional Navier-Stokes equations}
\Rjournal{Indiana Univ. Math. J.}
\Rvolume{57}
\Ryear{2008}
\Rnumber{6}
\Rpages{2643-2661}

\refer{Paper}{CW}
\Rauthor{Cao, Chongsheng and Jiahong Wu}
\Rtitle{Two regularity criteria for the 3D MHD equation}
\Rjournal{J. Differential Equations}
\Rvolume{248}
\Ryear{2010}
\Rnumber{9}
\Rpages{2263-2274}

\refer{Paper}{CZ}
\Rauthor{Chemin, Jean-Yves; Zhang, Ping}
\Rtitle{On the critical one component regularity for the 3-D Navier-Stokes equations}
\Rjournal{Ann. Sci. Éc. Norm. Supér.}
\Rvolume{49}
\Ryear{2016}
\Rnumber{1}
\Rpages{133-169}

\refer{Paper}{FO}
\Rauthor{Fan, Jishan; Ozawa, Tohru}
\Rtitle{Regularity criterion for weak solutions to the Navier-Stokes equations in terms of the gradient of the pressure}
\Rjournal{J. Inequal. Appl.}
\Rvolume{Art. ID 412678}
\Ryear{2008}
\Rpages{6}

\refer{Paper}{FJNZ}
\Rauthor{Fan, Jishan; Jiang, Song; Nakamura, Gen; Zhou, Yong}
\Rtitle{Logarithmically improved regularity criteria for the Navier-Stokes and MHD equations}
\Rjournal{J. Math. Fluid Mech.}
\Rvolume{13}
\Ryear{2011}
\Rnumber{4}
\Rpages{557–571}

\refer{Paper}{FQ1}
\Rauthor{Daoyuan Fang and  Chenyin Qian}
\Rtitle{Regularity criterion for 3D Navier-Stokes equations in Besov spaces}
\Rjournal{Commun. Pure Appl. Anal.}
\Rvolume{13}
\Ryear{2014}
\Rnumber{2}
\Rpages{585-603}

\refer{Paper}{FQ2}
\Rauthor{Fang, Daoyuan and Qian, Chenyin}
\Rtitle{The regularity criterion for 3{D} {N}avier-{S}tokes equations
              involving one velocity gradient component}
\Rjournal{Nonlinear Anal.
              International Multidisciplinary Journal}
\Rvolume{78}
\Ryear{2013}
\Rpages{86-103}

\refer{Paper}{G08}
\Rauthor{Sadek Gala}
\Rtitle{Regularity criterion on weak solutions to the Navier-Stokes equations}
\Rjournal{J Korean Math. Soc.}
\Rvolume{45}
\Ryear{2008}
\Rnumber{2}
\Rpages{537-558}

\refer{Paper}{GL}
\Rauthor{Lemari\'{e}-Rieusset, P. G. and Gala, S.}
\Rtitle{Multipliers between {S}obolev spaces and fractional
              differentiation}
\Rjournal{J. Math. Anal. Appl.}
\Rvolume{322}
\Ryear{2006}
\Rnumber{2}
\Rpages{1030-1054}

\refer{Paper}{GG}
\Rauthor{ZhengguangGuo and Sadek Gala}
\Rtitle{Remarks on logarithmical regularity criteria for the Navier-Stokes equations}
\Rjournal{J. Math. Phys.}
\Rvolume{52}
\Ryear{2011}
\Rnumber{6}
\Rpages{063503}

\refer{Paper}{He}
\Rauthor{He, Cheng}
\Rtitle{Regularity for solutions to the Navier-Stokes equations with one velocity component regular}
\Rjournal{Electron. J. Differential Equations 49}
\Rvolume{49}
\Ryear{2002}
\Rpages{1-13}
\Rnumber{29} 

\refer{Paper}{HG}
\Rauthor{Xiaowei He and Sadek Gala}
\Rtitle{Regularity criterion for weak solutions to the Navier-Stokes equations in terms of the pressure in the class $L^2(0, T; B_{\infty, \infty}^{-1}(\mathbb R^3))$}
\Rjournal{Nonlinear Anal. Real World Appl.
              Multidisciplinary Journal}
\Rvolume{12}
\Ryear{2011}
\Rnumber{6}
\Rpages{3602-3607}

\refer{Paper}{He1}
\Rauthor{He, Cheng}
\Rtitle{New sufficient conditions for regularity of solutions to the Navier-Stokes equations}
\Rjournal{Adv. Math. Sci. Appl.}
\Rvolume{12}
\Ryear{2002}
\Rnumber{02}
\Rpages{535-548}

\refer{Paper}{JZ1}
\Rauthor{Jia, Xuanji and Zhou, Yong}
\Rtitle{Remarks on regularity criteria for the {N}avier-{S}tokes
              equations via one velocity component}
\Rjournal{Nonlinear Anal. Real World Appl.
              Multidisciplinary Journal}
\Rvolume{15}
\Ryear{2014}
\Rpages{239-245}

\refer{Paper}{KOT}
\Rauthor{Hideo Kozono, Takayoshi Ogawa, Yasushi Taniuchi}
\Rtitle{The critical Sobolev inequalities in Besov spaces and regularity criterion to some semilinear evolution equations}
\Rjournal{Math. Z.}
\Rvolume{242}
\Ryear{2002}
\Rpages{251-278}

\refer{Paper}{Maz1}
\Rauthor{V. G. Maz'ya and T. O. Shaposhnikova}
\Rtitle{Theory of multipliers in spaces of differentiable functions,}
\Rjournal{Monographs and Studies in Mathematics, 23. Pitman (Advanced Publishing Program), Boston, MA,}
\Rvolume{}
\Ryear{1985}
\Rpages{xiii+344 pp}

\refer{Paper}{Maz2}
\Rauthor{V. G. Maz'ya}
\Rtitle{On the theory of the n-dimensional Schrödinger operator}
\Rjournal{Izv. Akad. Nauk SSSR (Ser. Mat.)}
\Rvolume{28}
\Ryear{1964}
\Rnumber{5}
\Rpages{1145-1172}

\refer{Paper}{PP}
\Rauthor{P. Penel and M. Pokorný}
\Rtitle{On anisotropic regularity criteria for the solutions to 3D Navier-Stokes equations}
\Rjournal{J. Math. Fluid Mech.}
\Rvolume{13}
\Ryear{2011}
\Rpages{341-353}

\refer{Book}{Pok}
\Rauthor{Milan Pokorný}
\Rtitle{On the result of He concerning the smoothness of solutions to the Navier-Stokes equations}
\Rjournal{Electron. J. Differ. Equ.}
\Rvolume{11} 
\Ryear{2003}
\Rpages{1-8}

\refer{Paper}{Ser}
\Rauthor{James Serrin}
\Rtitle{On the interior regularity of weak solutions of the Navier-Stokes equations}
\Rjournal{Arch. Ration. Mech. Anal.}
\Rvolume{9}
\Ryear{1962}
\Rpages{187-195}

\refer{Paper}{Ye1}
\Rauthor{Ye, Zhuan}
\Rtitle{Remarks on the regularity criterion to the 3{D}
              {N}avier-{S}tokes equations via one velocity component}
\Rjournal{Differential Integral Equations
              for Theory \& Applications}
\Rvolume{29}
\Ryear{2016}
\Rnumber{9-10}
\Rpages{957-976}

\refer{Paper}{Ye2}
\Rauthor{Ye, Zhuan}
\Rtitle{Remarks on the regularity criterion to the {N}avier-{S}tokes
              equations via the gradient of one velocity component}
\Rjournal{J. Math. Anal. Appl.}
\Rvolume{435}
\Ryear{2016}
\Rnumber{2}
\Rpages{1623-1633}

\refer{Paper}{ZC}
\Rauthor{Zhifei, Zhang; Qionglei, Chen}
\Rtitle{Regularity criterion via two components of vorticity on weak solutions to the Navier-Stokes equations in $\mathbb R^3$}
\Rjournal{J. Differential Equations}
\Rvolume{216}
\Ryear{2005}
\Rnumber{2}
\Rpages{470-481}

\refer{Paper}{ZY}
\Rauthor{Zujin Zhang and Xian Yang}
\Rtitle{A note on the regularity criterion for the 3D Navier-Stokes equations via the gradient of one velocity component}
\Rjournal{J. Math. Anal. Appl.}
\Rvolume{432}
\Ryear{2015}
\Rnumber{1}
\Rpages{603-611}

\refer{Paper}{Zhe}
\Rauthor{Zheng, Xiaoxin}
\Rtitle{A regularity criterion for the tridimensional
              {N}avier-{S}tokes equations in term of one velocity component}
\Rjournal{J. Differential Equations}
\Rvolume{256}
\Ryear{2014}
\Rnumber{1}
\Rpages{283-309}

\refer{Paper}{Z1}
\Rauthor{Zhou, Yong}
\Rtitle{A new regularity criterion for weak solutions to the Navier-Stokes equations}
\Rjournal{J. Math. Pures Appl.}
\Rvolume{9}
\Ryear{2005}
\Rnumber{84}
\Rpages{1496-1514}

\refer{Paper}{Z2}
\Rauthor {Zhou, Yong}
\Rtitle{A new regularity criterion for the Navier-Stokes equations in terms of the gradient of one velocity component}
\Rjournal{ Methods Appl. Anal.}
\Rvolume{9}
\Ryear{2002}
\Rnumber{4}
\Rpages{563-578}

\refer{Paper}{Z3}
\Rauthor{Zhou, Yong}
\Rtitle{Regularity criteria in terms of pressure for the 3-{D}
              {N}avier-{S}tokes equations in a generic domain}
\Rjournal{Math. Ann.}
\Rvolume{328}
\Ryear{2004}
\Rnumber{1-2}
\Rpages{173-192}

\refer{Paper}{Z4}
\Rauthor{Yong, Zhou}
\Rtitle{On regularity criteria in terms of pressure for the Navier-Stokes equations in $\mathbb R^3$}
\Rjournal{Proc. Amer. Math. Soc.}
\Rvolume{134}
\Ryear{2005}
\Rnumber{1}
\Rpages{149-156}

\refer{Paper}{Z5}
\Rauthor{Yong, Zhou}
\Rtitle{On a regularity criterion in terms of the gradient of pressure for the Navier-Stokes equations in $\mathbb R^N$}
\Rjournal{Z. Angew. Math. Phys.}
\Rvolume{57}
\Ryear{2006}
\Rpages{384-392}

\refer{Paper}{ZG1}
\Rauthor{Yong, Zhou and Sadek, Gala}
\Rtitle{Logarithmically improved regularity criteria for the Navier-Stokes equations in multiplier spaces}
\Rjournal{J. Math. Anal. Appl.}
\Rvolume{356}
\Ryear{2009}
\Rnumber{2}
\Rpages{498-501}

\refer{Paper}{ZP2}
\Rauthor{Zhou, Yong and Pokorn\'{y}, Milan}
\Rtitle{On a regularity criterion for the {N}avier-{S}tokes equations
              involving gradient of one velocity component}
\Rjournal{J. Math. Phys.}
\Rvolume{50}
\Ryear{2009}
\Rnumber{12}
\Rpages{123514, 11}

\refer{Paper}{ZP1}
\Rauthor{Zhou, Yong and Pokorn\'{y}, Milan}
\Rtitle{On the regularity of the solutions of the Navier-Stokes equations via one velocity component}
\Rjournal{Nonlinearity}
\Rvolume{23}
\Ryear{2010}
\Rnumber{5}
\Rpages{1097–1107}

\end{references}
